\documentclass[10pt]{article}
\usepackage{latexsym,amsfonts,amssymb}
\usepackage[dvips]{color}
\usepackage{colordvi,multicol}
\usepackage{amsmath}

\newtheorem{thm}{Theorem}[section]

\newtheorem{rem}[thm]{Remark}

\date{}
\begin{document}

\title{\bf A Note on the Gaussian Minimum Conjecture}
 \author{ Yang-Fan Zhong, Ting Ma,  Ze-Chun Hu\thanks{Corresponding author: zchu@scu.edu.cn}\\ \\
  {\small College of Mathematics, Sichuan  University,  China}}

\maketitle

\begin{abstract}
Let $n\geq 2$ and $(X_i,1\leq i\leq n)$ be a centered Gaussian random vector. The Gaussian minimum conjecture says that $E\left(\min_{1\leq i\leq n}|X_i|\right)\geq E\left(\min_{1\leq i\leq n}|Y_i|\right)$, where $Y_1,\ldots,Y_n$ are independent centered Gaussian random variables with $E(X_i^2)=E(Y_i^2)$ for any $i=1,\ldots,n$.
In this note, we will show that this conjecture holds if and only if  $n=2$.
\end{abstract}

\noindent  {\it MSC:} Primary 60E15; Secondary 62H12

\noindent  {\it Keywords:} Gaussian Minimum Conjecture; Gaussian random vector

\section{Introduction and main result}\setcounter{equation}{0}

Let $n\geq 2$ and $(X_i,1\leq i\leq n)$ be a centered Gaussian random vector. The well-known \u{S}id\'{a}k's inequality (\cite{Si67}, \cite{GL89}) says that
\begin{eqnarray}\label{1.1}
\mathbb{E}\left(\max_{1\leq i\leq n}|X_i|\right)\leq \mathbb{E}\left(\max_{1\leq i\leq n}|Y_i|\right),
\end{eqnarray}
where $Y_1,\ldots,Y_n$ are independent centered Gaussian random variables with $E(X_i^2)=E(Y_i^2)$ for any $i=1,\ldots,n$.

If we replace `max' by `min' in \u{S}id\'{a}k's inequality, Gordon et al. proved among other things that (\cite{GLSW05}, \cite{GLSW06})
\begin{eqnarray}\label{1.2}
\mathbb{E}\left(\min_{1\leq i\leq n}|X_i|\right)\geq \frac{1}{2}\mathbb{E}\left(\min_{1\leq i\leq n}|Y_i|\right).
\end{eqnarray}
Note that in \cite{GLSW05} and  \cite{GLSW06}, the authors proved the inequality \eqref{1.2} without any condition on the joint distribution of $(X_1,\ldots,X_n)$.

Wenbo V. Li and Qiman Shao conjectured that when $(X_i,1\leq i\leq n)$ is a centered Gaussian random vector, $\frac{1}{2}$ in the inequality \eqref{1.2} can be removed, i.e.
\begin{eqnarray}\label{main-ineq}
\mathbb{E}\left(\min_{1\leq i\leq n}|X_i|\right)\geq \mathbb{E}\left(\min_{1\leq i\leq n}|Y_i|\right),
\end{eqnarray}
which is called the {\it Gaussian Minimum Conjecture} (see \cite{Li12}, \cite{Sh16}).

Now we state the main result of this note.

\begin{thm}\label{thm-1}
The Gaussian Minimum Conjecture holds if and only if $n=2$.
\end{thm}

\begin{rem}
 Professor Qiman Shao told the third author that the Gaussian Minimum Conjecture for $n=2$ can be easily proved based on the following fact:
\begin{eqnarray*}
P(\min(|X_1|,|X_2|)>x)&=&1-P(|X_1|\leq x)-P(|X_2|\leq x)+P(|X_1|\leq x,|X_2|\leq x)\\
&\geq &1-P(|X_1|\leq x)-P(|X_2|\leq x)+P(|X_1|\leq x)P(|X_2|\leq x)\\
&=&1-P(|Y_1|\leq x)-P(|Y_2|\leq x)+P(|Y_1|\leq x)P(|Y_2|\leq x)\\
&=&1-P(|Y_1|\leq x)-P(|Y_2|\leq x)+P(|Y_1|\leq x, |Y_2|\leq x)\\
&=&P(\min(|Y_1|,|Y_2|)>x),
\end{eqnarray*}
where the Gaussian correlation inequality was used. In next section, we will give a different proof for $n=2$ based on an explicit formula for $E(\min(|X_1|,|X_2|))$ (see (\ref{2st-goal}) below).
\end{rem}

As to the minimum value point of $\mathbb{E}\left(\min_{1\leq i\leq n}|X_i|\right)$ with $E(X_i^2)=1,\forall i=1,\ldots,n$, for $n\geq 3$, we refer to \cite[Conjecture 5.1]{Litvak}. In addition, we would like to mention a related conjecture, the Gaussian random vector maximum (GRVM) conjecture, which asserts that among all centered Gaussian random vectors $X=(X_1,\ldots,X_n)$ with $E[X_i^2]=1,1\leq i\leq n$, the expectation $E[\max(X_1,\ldots,X_n)]$ is maximal if and only if all off-diagonal elements of the covariance matrix equal $-\frac{1}{n-1}$, where $n\geq 3$. The paper \cite{SHL20} proved the GRVM conjecture for $n=3,4$.

\section{Proof of Theorem \ref{thm-1}}\setcounter{equation}{0}

\subsection{Necessity}

 In this part, we will show that if $n\geq 3$, the  inequality (\ref{main-ineq}) does not hold for some $n$-dimensional centered Gaussian random vector $(X_1,\ldots,X_n)$ and independent centered Gaussian random variables $Y_1,\ldots,Y_n$ with $E(X_i^2)=E(Y_i^2),\forall i=1,\ldots,n$.

(i) $n=3$. R. van Handel gave the following counterexample for $n=3$ (see \cite{MZ11}). Let $Y_1,Y_2,Y_3$ be three independent standard Gaussian random variables. Define $X_1=(Y_1-Y_2)/\sqrt{2}, X_2=(Y_2-Y_3)/\sqrt{2},X_3=(Y_3-Y_1)/\sqrt{2}$. R. van Handel checked numerically that $0.32\sim \mathbb{E}\left(\min_{1\leq i\leq 3}|X_i|\right)<\mathbb{E}\left(\min_{1\leq i\leq 3}|Y_i|\right) \sim 0.34$, which tell us that (\ref{main-ineq}) does not hold for $n=3$. In fact, we can give the exact values of $\mathbb{E}\left(\min_{1\leq i\leq 3}|X_i|\right)$ and $\mathbb{E}\left(\min_{1\leq i\leq 3}|Y_i|\right)$.

At first, we calculate $\mathbb{E}\left(\min_{1\leq i\leq 3}|Y_i|\right)$. The density function $p_1(x,y,z)$ of $(Y_1,Y_2,Y_3)$ can be expressed by
$$
p_1(x,y,z)=\frac{1}{(2\pi)^{3/2}}e^{-\frac{x^2+y^2+z^2}{2}}.
$$
By the symmetry, we have
\begin{eqnarray}\label{Y-a}
\mathbb{E}\left(\min_{1\leq i\leq 3}|Y_i|\right)&=&
\int_{-\infty}^{\infty}\int_{-\infty}^{\infty}\int_{-\infty}^{\infty}(|x|\wedge |y|\wedge|z|)\frac{1}{(2\pi)^{3/2}}e^{-\frac{x^2+y^2+z^2}{2}}dxdydz\nonumber\\
&=&\frac{8}{(2\pi)^{3/2}}\int_0^{\infty}\int_0^{\infty}\int_0^{\infty}(x\wedge y\wedge z)e^{-\frac{x^2+y^2+z^2}{2}}dxdydz\nonumber\\
&=&\frac{8\cdot 3! }{(2\pi)^{3/2}}\int_0^{\infty}\left(\int_0^ze^{-\frac{y^2+z^2}{2}}\left(\int_0^y xe^{-\frac{x^2}{2}}dx\right)dy\right)dz\nonumber\\
&=&\frac{48}{(2\pi)^{3/2}}\int_0^{\infty}\left(\int_0^ze^{-\frac{y^2+z^2}{2}}
\left(1-e^{-\frac{y^2}{2}}\right)dy\right)dz\nonumber\\
&=&\frac{12\sqrt{2} }{\pi^{3/2}}\left[\int_0^{\infty}\left(\int_0^ze^{-\frac{y^2+z^2}{2}}dy\right)dz-
\int_0^{\infty}\left(\int_0^ze^{-\frac{2y^2+z^2}{2}}dy\right)dz\right].
\end{eqnarray}

Define a function
$$
F(a):=\int_0^{\infty}\left(\int_0^ze^{-\frac{ay^2+z^2}{2}}dy\right)dz,\ a>0,
$$
and a set
$$
D:=\{(y,z)\in \mathbf{R}^2: 0\leq y\leq z\}.
$$
Define a transformation
$$
(y,z):=T(u,v)=(u/\sqrt{a},v).
$$
Denote by $D_T$ the original image of $D$ under $T$. Then we have
\begin{eqnarray*}
D_T=\{(u,v)\in \mathbf{R}^2: T(u,v)\in D\}=\{(u,v)\in \mathbf{R}^2: 0\leq u/\sqrt{a}\leq v\}.
\end{eqnarray*}
Now,  we have
\begin{eqnarray}\label{Y-b}
F(a)&=&\iint_D e^{-\frac{ay^2+z^2}{2}}dydz=\frac{1}{\sqrt{a}}\iint_{D_T} e^{-\frac{u^2+v^2}{2}}dudv\nonumber\\
&=&\frac{1}{\sqrt{a}}\cdot \left(\frac{\pi}{2}-\arctan \frac{1}{\sqrt{a}}\right)\cdot \int_0^{\infty}e^{-\frac{r^2}{2}}rdr\nonumber\\
&=&\frac{1}{\sqrt{a}}\arctan\sqrt{a}.
\end{eqnarray}
By \eqref{Y-a} and \eqref{Y-b}, we get
\begin{eqnarray}\label{Y-c}
\mathbb{E}\left(\min_{1\leq i\leq 3}|Y_i|\right)
&=&\frac{12\sqrt{2} }{\pi^{3/2}}\left(F(1)-F(2)\right)\nonumber\\
&=&\frac{12\sqrt{2} }{\pi^{3/2}}\left(\arctan 1-\frac{1}{\sqrt{2}}\arctan\sqrt{2}\right)\nonumber\\
&=&\frac{12}{\pi^{3/2}}\left(\frac{\sqrt{2}\pi}{4}-\arctan\sqrt{2}\right).
\end{eqnarray}

Next, we come to calculate $\mathbb{E}\left(\min_{1\leq i\leq 3}|X_i|\right)$. Note that
$X_3=-(X_1+X_2)$, and the covariance matrix of $(X_1,X_2)$ is
$
\Sigma=\left(
\begin{array}{cc}
1 & -\frac{1}{2}\\
-\frac{1}{2}& 1
\end{array}
\right).
$
It follows that $\det\Sigma=\frac{3}{4}$, and
$
\Sigma^{-1}=\left(
\begin{array}{cc}
4/3 & 2/3\\
2/3& 4/3
\end{array}
\right).
$
Then  the density function $p_2(x,y)$ of $(X_1,X_2)$ can be expressed by
\begin{eqnarray*}
p_2(x,y)&=&\frac{1}{2\pi\sqrt{3/4}}e^{-\frac{\frac{4}{3}(x^2+y^2+xy)}{2}}
=\frac{1}{\sqrt{3}\pi}e^{-\frac{2(x^2+y^2+xy)}{3}}.
\end{eqnarray*}
By the symmetry, we have
\begin{eqnarray}\label{X-a}
\mathbb{E}\left(\min_{1\leq i\leq 3}|X_i|\right)&=&\frac{1}{\sqrt{3}\pi}\int_{-\infty}^{\infty}\int_{-\infty}^{\infty}
\left(|x|\wedge |y|\wedge |x+y|\right)e^{-\frac{2(x^2+y^2+xy)}{3}}dxdy\nonumber\\
&=&\frac{1}{\sqrt{3}\pi}\left[\int_0^{\infty}\int_0^{\infty}(x\wedge y)e^{-\frac{2(x^2+y^2+xy)}{3}}dxdy\right.
+\int_{-\infty}^0\int_{-\infty}^0\left((-x)\wedge (-y)\right) e^{-\frac{2(x^2+y^2+xy)}{3}}dxdy\nonumber\\
&&\quad\quad\quad\quad+\int_{-\infty}^0\int_0^{\infty}\left((-x)\wedge y\wedge |x+y|\right) e^{-\frac{2(x^2+y^2+xy)}{3}}dxdy\nonumber\\
&&\quad\quad\quad\quad\left.+\int_0^{\infty}\int_{-\infty}^0\left(x\wedge (-y)\wedge |x+y|\right) e^{-\frac{2(x^2+y^2+xy)}{3}}dxdy\right]\nonumber\\
&=&\frac{2}{\sqrt{3}\pi}\left[\int_0^{\infty}\int_0^{\infty}(x\wedge y)e^{-\frac{2(x^2+y^2+xy)}{3}}dxdy+\int_0^{\infty}\int_0^{\infty}\left(x\wedge y\wedge |x-y|\right) e^{-\frac{2(x^2+y^2-xy)}{3}}dxdy\right]\nonumber\\
&=&\frac{4}{\sqrt{3}\pi}\left[\int_0^{\infty}\left(\int_0^y xe^{-\frac{2(x^2+y^2+xy)}{3}}dx\right)dy+\int_0^{\infty}\left(\int_0^y\left(x\wedge (y-x)\right) e^{-\frac{2(x^2+y^2-xy)}{3}}dx\right)dy\right]\nonumber\\
&=&\frac{4}{\sqrt{3}\pi}\left[\int_0^{\infty}\left(\int_0^y xe^{-\frac{2(x^2+y^2+xy)}{3}}dx\right)dy\right.+\int_0^{\infty}\left(\int_0^{y/2}x e^{-\frac{2(x^2+y^2-xy)}{3}}dx\right)dy\nonumber\\
&&\quad\quad\quad\quad\left.+\int_0^{\infty}\left(\int_{y/2}^y (y-x)e^{-\frac{2(x^2+y^2-xy)}{3}}dx\right)dy\right]\nonumber\\
&=&\frac{4}{\sqrt{3}\pi}\left[\int_0^{\infty}\left(\int_0^y xe^{-\frac{2(x^2+y^2+xy)}{3}}dx\right)dy\right.+\int_0^{\infty}\left(\int_0^{y/2}x e^{-\frac{2(x^2+y^2-xy)}{3}}dx\right)dy\nonumber\\
&&\quad\quad\quad\quad+\int_0^{\infty}\left(\int_{y/2}^yy e^{-\frac{2(x^2+y^2-xy)}{3}}dx\right)dy\left.
-\int_0^{\infty}\left(\int_{y/2}^yxe^{-\frac{2(x^2+y^2-xy)}{3}}dx\right)dy\right]\nonumber\\
&=:&\frac{4}{\sqrt{3}\pi}\left[G_1+G_2+G_3-G_4\right],
\end{eqnarray}
where
\begin{eqnarray*}
&&G_1:=\int_0^{\infty}\left(\int_0^y xe^{-\frac{2(x^2+y^2+xy)}{3}}dx\right)dy,\quad\quad G_2:=\int_0^{\infty}\left(\int_0^{y/2}x e^{-\frac{2(x^2+y^2-xy)}{3}}dx\right)dy,\\
&&G_3:=\int_0^{\infty}\left(\int_{y/2}^yy e^{-\frac{2(x^2+y^2-xy)}{3}}dx\right)dy,\quad\quad G_4:=\int_0^{\infty}\left(\int_{y/2}^yxe^{-\frac{2(x^2+y^2-xy)}{3}}dx\right)dy.
\end{eqnarray*}

We have
\begin{eqnarray*}
G_1&=&\int_0^{\infty}\left(\int_0^y \left(x+\frac{y}{2}\right)e^{-\frac{2(x^2+y^2+xy)}{3}}dx\right)dy-
\frac{1}{2}\int_0^{\infty}\left(\int_0^y \left(y+\frac{x}{2}\right)e^{-\frac{2(x^2+y^2+xy)}{3}}dx\right)dy\\
&&+\frac{1}{4}\int_0^{\infty}\left(\int_0^y xe^{-\frac{2(x^2+y^2+xy)}{3}}dx\right)dy\\
&=&\frac{G_1}{4}+\int_0^{\infty}\left(\int_0^y \left(x+\frac{y}{2}\right)e^{-\frac{2\left(x+\frac{y}{2}\right)^2+\frac{3}{2}y^2}{3}}dx\right)dy-
\frac{1}{2}\int_0^{\infty}\left(\int_x^{\infty} \left(y+\frac{x}{2}\right)e^{-\frac{2\left(y+\frac{x}{2}\right)^2+\frac{3}{2}x^2}{3}}dy\right)dx\\
&=&\frac{G_1}{4}+\int_0^{\infty}e^{-\frac{y^2}{2}}\left(\int_{y/2}^{3y/2} ue^{-\frac{2u^2}{3}}du\right)dy-
\frac{1}{2}\int_0^{\infty}e^{-\frac{x^2}{2}}\left(\int_{3x/2}^{\infty} ve^{-\frac{2v^2}{3}}dv\right)dx\\
&=&\frac{G_1}{4}+\frac{3}{4}\left[\int_0^{\infty}e^{-\frac{y^2}{2}}
\left(e^{-\frac{y^2}{6}}
-e^{-\frac{3y^2}{2}}\right)dy-
\frac{1}{2}\int_0^{\infty}e^{-\frac{x^2}{2}}\cdot e^{-\frac{3x^2}{2}}dx\right]\\
&=&\frac{G_1}{4}+\frac{3}{8}\left(2\int_0^{\infty}e^{-\frac{2y^2}{3}}dy
-3\int_0^{\infty}e^{-2y^2}dy\right).
\end{eqnarray*}
It follows that
\begin{eqnarray}\label{G-1}
G_1=\frac{1}{2}\left(2\int_0^{\infty}e^{-\frac{2y^2}{3}}dy
-3\int_0^{\infty}e^{-2y^2}dy\right).
\end{eqnarray}

We have
\begin{eqnarray*}
G_2&=&\int_0^{\infty}\left(\int_0^{y/2}\left(x-\frac{y}{2}\right) e^{-\frac{2(x^2+y^2-xy)}{3}}dx\right)dy+\frac{1}{2}\int_0^{\infty}\left(\int_0^{y/2}
\left(y-\frac{x}{2}\right) e^{-\frac{2(x^2+y^2-xy)}{3}}dx\right)dy\\
&&+\frac{1}{4}\int_0^{\infty}\left(\int_0^{y/2}x e^{-\frac{2(x^2+y^2-xy)}{3}}dx\right)dy\\
&=&\frac{G_2}{4}+\int_0^{\infty}\left(\int_0^{y/2}\left(x-\frac{y}{2}\right) e^{-\frac{2\left(x-\frac{y}{2}\right)^2+\frac{3y^2}{2}}{3}}dx\right)dy
+\frac{1}{2}\int_0^{\infty}\left(\int_{2x}^{\infty}
\left(y-\frac{x}{2}\right) e^{-\frac{2\left(y-\frac{x}{2}\right)^2+\frac{3x^2}{2}}{3}}dy\right)dx\\
&=&\frac{G_2}{4}+\int_0^{\infty}e^{-\frac{y^2}{2}}
\left(\int_{-y/2}^0u e^{-\frac{2u^2}{3}}du\right)dy
+\frac{1}{2}\int_0^{\infty}e^{-\frac{x^2}{2}}\left(\int_{3x/2}^{\infty}
v e^{-\frac{2v^2}{3}}dv\right)dx\\
&=&\frac{G_2}{4}+\frac{3}{4}\left[\int_0^{\infty}e^{-\frac{y^2}{2}}
\left(e^{-\frac{y^2}{6}}-1\right)dy
+\frac{1}{2}\int_0^{\infty}e^{-\frac{x^2}{2}}\cdot e^{-\frac{3x^2}{2}}dx\right]\\
&=&\frac{G_2}{4}+\frac{3}{8}\left(2\int_0^{\infty}e^{-\frac{2y^2}{3}}dy
-2\int_0^{\infty}e^{-\frac{y^2}{2}}dy+\int_0^{\infty}e^{-2x^2}dx\right).
\end{eqnarray*}
It follows that
\begin{eqnarray}\label{G-2}
G_2&=&\frac{1}{2}\left(2\int_0^{\infty}e^{-\frac{2y^2}{3}}dy
-2\int_0^{\infty}e^{-\frac{y^2}{2}}dy+\int_0^{\infty}e^{-2x^2}dx\right).
\end{eqnarray}

We have
\begin{eqnarray*}
G_3&=&\int_0^{\infty}\left(\int_{y/2}^y\left(y-\frac{x}{2}\right) e^{-\frac{2(x^2+y^2-xy)}{3}}dx\right)dy+
\frac{1}{2}\int_0^{\infty}\left(\int_{y/2}^y\left(x-\frac{y}{2}\right) e^{-\frac{2(x^2+y^2-xy)}{3}}dx\right)dy\\
&&+\frac{1}{4}\int_0^{\infty}\left(\int_{y/2}^yy e^{-\frac{2(x^2+y^2-xy)}{3}}dx\right)dy\\
&=&\frac{G_3}{4}+\int_0^{\infty}\left(\int_{x}^{2x}\left(y-\frac{x}{2}\right) e^{-\frac{2\left(y-\frac{x}{2}\right)^2+\frac{3x^2}{2}}{3}}dy\right)dx+
\frac{1}{2}\int_0^{\infty}\left(\int_{y/2}^y\left(x-\frac{y}{2}\right) e^{-\frac{2\left(x-\frac{y}{2}\right)^2+\frac{3y^2}{2}}{3}}dx\right)dy\\
&=&\frac{G_3}{4}+\int_0^{\infty}e^{-\frac{x^2}{2}}\left(\int_{x/2}^{3x/2}
u e^{-\frac{2u^2}{3}}du\right)dx+
\frac{1}{2}\int_0^{\infty}e^{-\frac{y^2}{2}}\left(\int_0^{y/2}v e^{-\frac{2v^2}{3}}dv\right)dy\\
&=&\frac{G_3}{4}+\frac{3}{4}\left[\int_0^{\infty}e^{-\frac{x^2}{2}}
\left(e^{-\frac{x^2}{6}}-e^{-\frac{3x^2}{2}}\right)dx+
\frac{1}{2}\int_0^{\infty}e^{-\frac{y^2}{2}}\left(1-e^{-\frac{y^2}{6}}\right)dy\right]\\
&=&\frac{G_3}{4}+\frac{3}{8}\left(\int_0^{\infty}e^{-\frac{y^2}{2}}dy+
\int_0^{\infty}e^{-\frac{2x^2}{3}}dx-2\int_0^{\infty}e^{-2x^2}dx\right).
\end{eqnarray*}
It follows that
\begin{eqnarray}\label{G-3}
G_3&=&\frac{1}{2}\left(\int_0^{\infty}e^{-\frac{y^2}{2}}dy+
\int_0^{\infty}e^{-\frac{2x^2}{3}}dx-2\int_0^{\infty}e^{-2x^2}dx\right).
\end{eqnarray}

We have
\begin{eqnarray*}
G_4&=&\int_0^{\infty}\left(\int_{y/2}^y\left(x-\frac{y}{2}\right)
e^{-\frac{2(x^2+y^2-xy)}{3}}dx\right)dy+
\frac{1}{2}\int_0^{\infty}\left(\int_{y/2}^y\left(y-\frac{x}{2}\right)
e^{-\frac{2(x^2+y^2-xy)}{3}}dx\right)dy\\
&&+\frac{1}{4}\int_0^{\infty}\left(\int_{y/2}^yxe^{-\frac{2(x^2+y^2-xy)}{3}}dx\right)dy\\
&=&\frac{G_4}{4}+\int_0^{\infty}\left(\int_{y/2}^y\left(x-\frac{y}{2}\right)
e^{-\frac{2\left(x-\frac{y}{2}\right)^2+\frac{3y^2}{2}}{3}}dx\right)dy+
\frac{1}{2}\int_0^{\infty}\left(\int_x^{2x}\left(y-\frac{x}{2}\right)
e^{-\frac{2\left(y-\frac{x}{2}\right)^2+\frac{3x^2}{2}}{3}}dy\right)dx\\
&=&\frac{G_4}{4}+\int_0^{\infty}e^{-\frac{y^2}{2}}
\left(\int_0^{y/2}ue^{-\frac{2u^2}{3}}du\right)dy+
\frac{1}{2}\int_0^{\infty}e^{-\frac{x^2}{2}}\left(\int_{x/2}^{3x/2}v
e^{-\frac{2v^2}{3}}dv\right)dx\\
&=&\frac{G_4}{4}+\frac{3}{4}\left[\int_0^{\infty}e^{-\frac{y^2}{2}}
\left(1-e^{-\frac{y^2}{6}}\right)dy+
\frac{1}{2}\int_0^{\infty}e^{-\frac{x^2}{2}}\left(e^{-\frac{x^2}{6}}
-e^{-\frac{3x^2}{2}}\right)dx\right]\\
&=&\frac{G_4}{4}+\frac{3}{8}\left(2\int_0^{\infty}e^{-\frac{y^2}{2}}dy
-\int_0^{\infty}e^{-\frac{2y^2}{3}}dy-\int_0^{\infty}e^{-2y^2}dy\right).
\end{eqnarray*}
It follows that
\begin{eqnarray}\label{G-4}
G_4&=&\frac{1}{2}\left(2\int_0^{\infty}e^{-\frac{y^2}{2}}dy
-\int_0^{\infty}e^{-\frac{2y^2}{3}}dy-\int_0^{\infty}e^{-2y^2}dy\right).
\end{eqnarray}

By \eqref{X-a}-\eqref{G-4}, we obtain
\begin{eqnarray}\label{XX}
\mathbb{E}\left(\min_{1\leq i\leq 3}|X_i|\right)
&=&\frac{6}{\sqrt{3}\pi}\left(2\int_0^{\infty}e^{-\frac{2y^2}{3}}dy
-\int_0^{\infty}e^{-2y^2}dy-\int_0^{\infty}e^{-\frac{y^2}{2}}dy\right)\nonumber\\
&=&\frac{6}{\sqrt{3}\pi}\cdot \frac{\sqrt{2\pi}}{2}\left(2\sqrt{\frac{3}{4}}-\sqrt{\frac{1}{4}}-1\right)\nonumber\\
&=&\frac{3(2-\sqrt{3})}{\sqrt{2\pi}}.
\end{eqnarray}

Hence, we get
\begin{eqnarray*}
\mathbb{E}\left(\min_{1\leq i\leq 3}|Y_i|\right)-\mathbb{E}\left(\min_{1\leq i\leq 3}|X_i|\right)&=&\frac{12}{\pi^{3/2}}
\left(\frac{\sqrt{2\pi}}{4}-\arctan\sqrt{2}\right)-\frac{3(2-\sqrt{3})}{\sqrt{2\pi}}\\
&=&\frac{12}{\pi^{3/2}}\left(\frac{\sqrt{6}\pi}{8}-\arctan\sqrt{2}\right)>0,
\end{eqnarray*}
since
\begin{eqnarray*}
\tan\left(\arctan\sqrt{2}\right)=\sqrt{2}\sim 1.414,\ \mbox{and}\ \tan \frac{\sqrt{6}\pi}{8}\sim 1.434.
\end{eqnarray*}

(ii) $n\geq 4$.   Without loss of generality, we only consider the case that $n=4$. We use proof by contradiction. Suppose that (\ref{main-ineq}) holds for $n=4$. Let $Y_i,X_i,i=1,2,3$ be the same as in the above example. Let $Y_4$ be a standard Gaussian random variable independent of $(Y_1,Y_2,Y_3)$. Then by the assumption, for any $a>0$, we have
\begin{eqnarray*}
\mathbb{E}\left((\min_{1\leq i\leq 3}|X_i|)\wedge |aY_4|\right)\geq \mathbb{E}\left((\min_{1\leq i\leq 3}|Y_i|)\wedge |aY_4|\right).
\end{eqnarray*}
Letting $a\to\infty$, by the monotone convergence theorem,  we obtain that
\begin{eqnarray*}
\mathbb{E}\left(\min_{1\leq i\leq 3}|X_i|\right)\geq \mathbb{E}\left(\min_{1\leq i\leq 3}|Y_i|\right).
\end{eqnarray*}
It is a contradiction. Hence for any $M>0$, there exists $a_0>M$ such that
\begin{eqnarray*}
\mathbb{E}\left((\min_{1\leq i\leq 3}|X_i|)\wedge |a_0Y_4|\right)< \mathbb{E}\left((\min_{1\leq i\leq 3}|Y_i|)\wedge |a_0Y_4|\right).
\end{eqnarray*}

\subsection{Sufficiency}

 In this part, we will show that the inequality (\ref{main-ineq}) holds if $n=2$.

Write $X_1=x_1f_1,X_2=x_2f_2$, where both $f_1$ and $f_2$ have the standard normal distribution $N(0,1)$.  Without loss of generality, we can assume that $x_1,x_2>0$. Further
we can assume that $x_1=1,x_2=a\in (0,1]$.

Denote by $\Sigma_1$ the covariance matrix of $(f_1,f_2)$. By the symmetry of the distribution $N(0,1)$, we can assume that
$
\Sigma_1=\left(
\begin{array}{cc}
1&\rho\\
\rho&1
\end{array}
\right),
$
where $0\leq \rho\leq 1$. Then the covariance matrix $\Sigma_2$ of $(f_1,af_2)$ can be expressed by
$
\Sigma_2=\left(
\begin{array}{cc}
1&a\rho\\
a\rho&a^2
\end{array}
\right).
$
It follows that if $\rho\in [0,1)$, then
\begin{eqnarray*}
\Sigma_2^{-1}=\frac{1}{a^2(1-\rho^2)}\left(
\begin{array}{cc}
a^2&-a\rho\\
-a\rho&1
\end{array}
\right),
\end{eqnarray*}
and thus the density function of $(f_1,af_2)$ is
$$
p(x,y)=\frac{1}{2\pi\sqrt{a^2(1-\rho^2)}}\exp\left(-\frac{a^2x^2+y^2-2a\rho xy}{2a^2(1-\rho^2)}\right).
$$

At first, we assume that $\rho\in [0,1)$. By the symmetry, we have
\begin{eqnarray*}
&&\mathbb{E}(|f_1|\wedge |af_2|)\\
&&=\int_{-\infty}^{\infty}\int_{-\infty}^{\infty}(|x|\wedge |y|)\frac{1}{2\pi\sqrt{a^2(1-\rho^2)}}e^{-\frac{a^2x^2+y^2-2a\rho xy}{2a^2(1-\rho^2)}}dxdy\nonumber\\
&&=\frac{1}{2\pi\sqrt{a^2(1-\rho^2)}}\left[\int_0^{\infty}\int_0^{\infty}(x\wedge y)e^{-\frac{a^2x^2+y^2-2a\rho xy}{2a^2(1-\rho^2)}}dxdy\right.\\
&&\quad+\int_{-\infty}^{0}\int_{-\infty}^0((-x)\wedge (-y))e^{-\frac{a^2x^2+y^2-2a\rho xy}{2a^2(1-\rho^2)}}dxdy+\int_0^{\infty}\int_{-\infty}^0(x\wedge (-y))e^{-\frac{a^2x^2+y^2-2a\rho xy}{2a^2(1-\rho^2)}}dxdy\\
&&\quad\left.+\int_{-\infty}^0\int_0^{\infty}((-x)\wedge y)e^{-\frac{a^2x^2+y^2-2a\rho xy}{2a^2(1-\rho^2)}}dxdy\right]\\
&&=\frac{1}{\pi\sqrt{a^2(1-\rho^2)}}\left[\int_0^{\infty}\int_0^{\infty}(x\wedge y)e^{-\frac{a^2x^2+y^2-2a\rho xy}{2a^2(1-\rho^2)}}dxdy+\int_0^{\infty}\int_0^{\infty}(x\wedge y)e^{-\frac{a^2x^2+y^2+2a\rho xy}{2a^2(1-\rho^2)}}dxdy\right].
\end{eqnarray*}
Define
$$
I(\vartheta):=\int_0^{\infty}\int_0^{\infty}(x\wedge y)e^{-\frac{a^2x^2+y^2+2a\vartheta xy}{2a^2(1-\rho^2)}}dxdy.
$$
Then
\begin{eqnarray*}
\mathbb{E}(|f_1|\wedge |af_2|)=\frac{1}{\pi\sqrt{a^2(1-\rho^2)}}\left[I(\rho)+I(-\rho)\right].
\end{eqnarray*}

We have
\begin{eqnarray*}
I(\vartheta)&=&\int_0^{\infty}\int_0^{\infty}(x\wedge y)e^{-\frac{a^2x^2+y^2+2a\vartheta xy}{2a^2(1-\rho^2)}}dxdy\\
&=&\int_0^{\infty}\left(\int_0^yxe^{-\frac{a^2x^2+y^2+2a\vartheta xy}{2a^2(1-\rho^2)}}dx\right)dy+\int_0^{\infty}\left(\int_0^xye^{-\frac{a^2x^2+y^2+2a\vartheta xy}{2a^2(1-\rho^2)}}dy\right)dx.
\end{eqnarray*}
Define
$$
J(\alpha,\beta,\gamma):=\int_0^{\infty}\left(\int_0^yxe^{-\frac{\alpha x^2+\beta y^2+2\gamma xy}{2a^2(1-\rho^2)}}dx\right)dy,
$$
where $\alpha>0,\beta>0,\alpha\beta-\gamma^2>0$. Then we have
$$
I(\vartheta)=J(a^2,1,a\vartheta)+J(1,a^2,a\vartheta),
$$
and thus
\begin{eqnarray}\label{1st-goal}
\mathbb{E}(|f_1|\wedge |af_2|)=\frac{1}{\pi\sqrt{a^2(1-\rho^2)}}\left[J(a^2,1,a\rho)+J(1,a^2,a\rho)
+J(a^2,1,-a\rho)+J(1,a^2,-a\rho)\right].\ \
\end{eqnarray}

In the following, we come to calculate the function $J(\alpha,\beta,\gamma)$. We have
\begin{eqnarray}\label{1st-goal-function}
J(\alpha,\beta,\gamma)&=&\int_0^{\infty}\left(\int_0^yxe^{-\frac{\alpha x^2+\beta y^2+2\gamma xy}{2a^2(1-\rho^2)}}dx\right)dy\nonumber\\
&=&\int_0^{\infty}\left(\int_0^y\left(x+\frac{\gamma}{\alpha}y\right)e^{-\frac{\alpha x^2+\beta y^2+2\gamma xy}{2a^2(1-\rho^2)}}dx\right)dy
-\frac{\gamma}{\alpha}\int_0^{\infty}\left(\int_0^yye^{-\frac{\alpha x^2+\beta y^2+2\gamma xy}{2a^2(1-\rho^2)}}dx\right)dy\nonumber\\
&=:&J_1(\alpha,\beta,\gamma)-J_2(\alpha,\beta,\gamma),
\end{eqnarray}
where
\begin{eqnarray*}
J_1(\alpha,\beta,\gamma)&:=&\int_0^{\infty}\left(\int_0^y\left(x+\frac{\gamma}{\alpha}y\right)e^{-\frac{\alpha x^2+\beta y^2+2\gamma xy}{2a^2(1-\rho^2)}}dx\right)dy,\\
J_2(\alpha,\beta,\gamma)&:=&\frac{\gamma}{\alpha}\int_0^{\infty}\left(\int_0^yye^{-\frac{\alpha x^2+\beta y^2+2\gamma xy}{2a^2(1-\rho^2)}}dx\right)dy.
\end{eqnarray*}
We have
\begin{eqnarray}\label{J-1}
J_1(\alpha,\beta,\gamma)&=&\int_0^{\infty}\left(\int_0^y\left(x+\frac{\gamma}{\alpha}y\right)
e^{-\frac{\alpha \left(x+\frac{\gamma}{\alpha}y\right)^2+\frac{\alpha\beta-\gamma^2}{\alpha} y^2}{2a^2(1-\rho^2)}}dx\right)dy\nonumber\\
&=&\int_0^{\infty}e^{-\frac{(\alpha\beta-\gamma^2)y^2}{2a^2(1-\rho^2)\alpha}}
\left(\int_{\frac{\gamma}{\alpha}y}^{\left(1+\frac{\gamma}{\alpha}\right)y}
ue^{-\frac{\alpha u^2}{2a^2(1-\rho^2)}}du\right)dy\nonumber\\
&=&\frac{a^2(1-\rho^2)}{\alpha}\int_0^{\infty}e^{-\frac{(\alpha\beta-\gamma^2)y^2}{2a^2(1-\rho^2)\alpha}}
\left(e^{-\frac{\gamma^2y^2}{2a^2(1-\rho^2)\alpha}}
-e^{-\frac{\left(\alpha+\gamma\right)^2y^2}{2a^2(1-\rho^2)\alpha}}\right)dy\nonumber\\
&=&\frac{a^2(1-\rho^2)}{\alpha}\int_0^{\infty}
\left(e^{-\frac{\beta y^2}{2a^2(1-\rho^2)}}
-e^{-\frac{\left(\alpha+\beta+2\gamma\right)y^2}{2a^2(1-\rho^2)}}\right)dy\nonumber\\
&=&\frac{a^2(1-\rho^2)}{\alpha}\cdot \frac{\sqrt{2\pi}}{2}\left(\frac{a\sqrt{1-\rho^2}}{\sqrt{\beta}}-
\frac{a\sqrt{1-\rho^2}}{\sqrt{\alpha+\beta+2\gamma}}\right)\nonumber\\
&=&\frac{\sqrt{2\pi}a^3(1-\rho^2)^{\frac{3}{2}}}{2\alpha}
\left(\frac{1}{\sqrt{\beta}}-
\frac{1}{\sqrt{\alpha+\beta+2\gamma}}\right),
\end{eqnarray}
and
\begin{eqnarray}\label{J-2}
J_2(\alpha,\beta,\gamma)
&=&\frac{\gamma}{\alpha}\int_0^{\infty}\left(\int_0^y\left(y+\frac{\gamma}{\beta}x\right)e^{-\frac{\alpha x^2+\beta y^2+2\gamma xy}{2a^2(1-\rho^2)}}dx\right)dy\nonumber\\
&&
-\frac{\gamma}{\alpha}\int_0^{\infty}\left(\int_0^y\left(\frac{\gamma}{\beta}x\right)e^{-\frac{\alpha x^2+\beta y^2+2\gamma xy}{2a^2(1-\rho^2)}}dx\right)dy\nonumber\\
&=&-\frac{\gamma^2}{\alpha\beta}\int_0^{\infty}\left(\int_0^yxe^{-\frac{\alpha x^2+\beta y^2+2\gamma xy}{2a^2(1-\rho^2)}}dx\right)dy\nonumber\\
&&+\frac{\gamma}{\alpha}\int_0^{\infty}\left(\int_x^{\infty}
\left(y+\frac{\gamma}{\beta}x\right)
e^{-\frac{\beta\left(y+\frac{\gamma}{\beta}x\right)^2+
\frac{\alpha\beta-\gamma^2}{\beta}x^2}{2a^2(1-\rho^2)}}dy\right)dx\nonumber\\
&=&-\frac{\gamma^2}{\alpha\beta}J(\alpha,\beta,\gamma)+
\frac{\gamma}{\alpha}\int_0^{\infty}
e^{-\frac{(\alpha\beta-\gamma^2)x^2}{2a^2(1-\rho^2)\beta}}
\left(\int_{\frac{(\beta+\gamma)x}{\beta}}^{\infty}
ue^{-\frac{\beta u^2}{2a^2(1-\rho^2)}}du\right)dx\nonumber\\
&=&-\frac{\gamma^2}{\alpha\beta}J(\alpha,\beta,\gamma)+
\frac{\gamma}{\alpha}\int_0^{\infty}
e^{-\frac{(\alpha\beta-\gamma^2)x^2}{2a^2(1-\rho^2)\beta}}\cdot \frac{a^2(1-\rho^2)}{\beta}e^{-\frac{(\beta+\gamma)^2x^2}{2a^2(1-\rho^2)\beta}}dx\nonumber\\
&=&-\frac{\gamma^2}{\alpha\beta}J(\alpha,\beta,\gamma)+
\frac{a^2(1-\rho^2)\gamma}{\alpha\beta}\int_0^{\infty}
e^{-\frac{(\alpha+\beta+2\gamma)x^2}{2a^2(1-\rho^2)}}dx\nonumber\\
&=&-\frac{\gamma^2}{\alpha\beta}J(\alpha,\beta,\gamma)+
\frac{\sqrt{2\pi}a^3(1-\rho^2)^{\frac{3}{2}}\gamma}{2\alpha\beta\sqrt{\alpha+\beta+2\gamma}}.
\end{eqnarray}
By \eqref{1st-goal-function} - \eqref{J-2}, we get
\begin{eqnarray*}
J(\alpha,\beta,\gamma)&=&\frac{1}{1-\frac{\gamma^2}{\alpha\beta}}\left[\frac{\sqrt{2\pi}a^3(1-\rho^2)^{\frac{3}{2}}}{2\alpha}
\left(\frac{1}{\sqrt{\beta}}-
\frac{1}{\sqrt{\alpha+\beta+2\gamma}}\right)-
\frac{\sqrt{2\pi}a^3(1-\rho^2)^{\frac{3}{2}}\gamma}{2\alpha\beta\sqrt{\alpha+\beta+2\gamma}}
\right]\\
&=&\frac{\sqrt{2\pi}a^3(1-\rho^2)^{\frac{3}{2}}}{2}\cdot
\frac{\beta}{\alpha\beta-\gamma^2}\left(\frac{1}{\sqrt{\beta}}-
\frac{\beta+\gamma}{\beta\sqrt{\alpha+\beta+2\gamma}}\right),
\end{eqnarray*}
which together with \eqref{1st-goal} implies that
\begin{eqnarray}\label{2st-goal}
&&\mathbb{E}(|f_1|\wedge |af_2|)\nonumber\\
&&=\frac{1}{\pi\sqrt{a^2(1-\rho^2)}}\left[J(a^2,1,a\rho)+J(1,a^2,a\rho)
+J(a^2,1,-a\rho)+J(1,a^2,-a\rho)\right]\nonumber\\
&&=\frac{1}{\pi\sqrt{a^2(1-\rho^2)}}\cdot \frac{\sqrt{2\pi}a^3(1-\rho^2)^{\frac{3}{2}}}{2}\left[
\frac{1}{a^2(1-\rho^2)}\left(1-
\frac{1+ap}{\sqrt{a^2+1+2a\rho}}\right)\right.\nonumber\\
&&\quad+\frac{a^2}{a^2(1-\rho^2)}\left(\frac{1}{a}-
\frac{a^2+a\rho}{a^2\sqrt{1+a^2+2a\rho}}\right)
+\frac{1}{a^2(1-\rho^2)}\left(1-
\frac{1-a\rho}{\sqrt{a^2+1-2a\rho}}\right)\nonumber\\
&&\quad\left.+\frac{a^2}{a^2(1-\rho^2)}\left(\frac{1}{a}-
\frac{a^2-a\rho}{a^2\sqrt{a^2+1-2a\rho}}\right)\right]\nonumber\\
&&=\frac{a^2(1-\rho^2)}{\sqrt{2\pi}}\cdot \frac{2(1+a)-\sqrt{1+a^2+2a\rho}-\sqrt{1+a^2-2a\rho}}{a^2(1-\rho^2)}\nonumber\\
&&=\frac{2(1+a)-\sqrt{1+a^2+2a\rho}-\sqrt{1+a^2-2a\rho}}{\sqrt{2\pi}}.
\end{eqnarray}

If $\rho=1$, then $f_2=f_1$ a.s. Note that $a\in(0,1]$. Then we have
\begin{eqnarray*}
\mathbb{E}(|f_1|\wedge |af_2|)&=&a\mathbb{E}(|f_1|)\\
&=&a\int_{-\infty}^{\infty}|x|\cdot \frac{1}{\sqrt{2\pi}}e^{-\frac{x^2}{2}}dx\\
&=&\frac{2a}{\sqrt{2\pi}}\int_0^{\infty}xe^{-\frac{x^2}{2}}dx\\
&=&\frac{2a}{\sqrt{2\pi}}.
\end{eqnarray*}
In addition, if $\rho=1$, we have
$$
\frac{2(1+a)-\sqrt{1+a^2+2a\rho}-\sqrt{1+a^2-2a\rho}}{\sqrt{2\pi}}=\frac{2a}{\sqrt{2\pi}}.
$$
Hence \eqref{2st-goal} holds for any $\rho\in [0,1]$.

For any $a\in (0,1]$, and any $\rho\in (0,1)$, we have
\begin{eqnarray*}
\frac{d\left[\frac{2(1+a)-\sqrt{1+a^2+2a\rho}-\sqrt{1+a^2-2a\rho}}{\sqrt{2\pi}}\right]}{d\rho}
=\frac{a}{\sqrt{2\pi}}\left(\frac{1}{\sqrt{1+a^2-2a\rho}}-
\frac{1}{\sqrt{1+a^2+2a\rho}}\right)> 0.
\end{eqnarray*}
Hence for any $a\in (0,1]$, $\frac{2(1+a)-\sqrt{1+a^2+2a\rho}-\sqrt{1+a^2-2a\rho}}{\sqrt{2\pi}}$ is a strictly increasing function in $\rho\in [0,1]$. Hence  it reaches its minimum value at $\rho=0$, i.e. the inequality (\ref{main-ineq}) holds. The proof is complete.

\bigskip

{ \noindent {\bf\large Acknowledgments}\ We thank Guolie Lan, Qiman Shao and Wei Sun for helpful discussions. This work was supported by the National Natural Science Foundation of China (No. 11771309).

\end{document}